\documentclass{amsart}
\usepackage{graphicx}
\newcommand{\C}{\mathbb{C}}
\newcommand{\Z}{\mathbb{Z}}

\input{amssym.def}
\input{amssym.tex}

\newcommand{\ben}{\begin{enumerate}}
\newcommand{\een}{\end{enumerate}}

\newcommand{\be}{\begin{enumerate}}
\newcommand{\ee}{\end{enumerate}}
\newcommand{\bq}{\begin{eqnarray*}}
\newcommand{\eq}{\end{eqnarray*}}

\newcommand{\disp}{\displaystyle}
\begin{document}
\title[THE  FUZZY SUBGROUPS RESULTS      ]{  THE  FUZZY SUBGROUPS   RESULTS  INVOLVING  MULTIPLE   SUMS    }
\author{  S. A.  Adebisi$^{1}$ ,    M.  Ogiugo $^{2} $ \&    M.   EniOluwafe$^{2}$  }
\address{$^{1}$Department of Mathematics , Faculty of Science, University of  Lagos,  Nigeria. Email :$\ adesinasunday@yahoo.com $,phone number : +2\;3\;4\;7\;0\;4\;1\;6\;3\;9\;0\;1\;3  \; \; \;\;\;\; \;\;\; \;\;\; \;\;\; \;\;\; \;\;\; \;\;\; \;\;\; \;\;\; \;\;\; \;\;\; \;\;\; \;\;\; \;\;\; \;\;\; \;\;\; \;\;\; \;\;\; \;\;\;\;\; \;\;\; \;\;\; \;\;\; \;\;\; \;\;\; \;\;\; \;\;\; \;\;\; 
$^{2}$Department of mathematics,Faculty of Science, university of Ibadan. Nigiria }
\keywords{Finite $p$-Groups, Nilpotent Group, Fuzzy subgroups, Dihedral Group, Inclusion-Exclusion Principle,Maximal subgroups. AMS Mathematics Subject Classification 2020: Primary: 08A72,  20D15, 60A86 .    Secondary :   20N25  \\  ORCID of the corresponding author: https://orcid.org/0000-0003-3766-0910}
\begin{abstract}     The theory of fuzzy sets has a wide range of applications, one of which is that of fuzzy groups . The fuzzy sets were introduced by  Zadeh. Even though, the story of fuzzy logic started much earlier, it was specially designed mathematically to represent uncertainty and vagueness. It was also, to provide formalized tools for dealing with the imprecision intrinsic to many problems. A group is said  to be  nilpotent if it has a normal series of a finite length $ n $. By this notion, every finite p-group is nilpotent. Nilpotent structures  such as the  p-groups, have normal series of finite length. Any finite p-group has many normal subgroups and consequently, the phenomenon of large number of non-isomorphic subgroups of a given order. This makes it an ideal object for combinatorial and cohomological investigations. Cartesian product (otherwise known as the product set) plays vital roles in the course of synthesizing the abstract groups. Previous studies have determined the number of distinct fuzzy subgroups of various finite p-groups including those of square-free order. However, not much work has been done on the fuzzy subgroup classification for the nilpotent groups formed from the Cartesian products of p-groups through their computations. This work is therefore designed to classify the nilpotent groups formed from the Cartesian products of $ p $-groups through their computations. In this paper, the Cartesian products of  $ p $-groups were taken to obtain nilpotent groups. the explicit formulae is given for the number of distinct fuzzy subgroups of the Cartesian product of the dihedral group of order eight with a cyclic group of order of an $n$ power of two for, which $n$ is not less than three . . \\
\end{abstract}
\maketitle
\title

 \section{  Introduction  }  

The aspect of pure Mathematics has undergone a lot of dynamic developments over the years .  Concerning   the  theory of  fuzzy group , the classification, most especially the finite $p$-groups cannot be overlooked.   For instance, many researchers have treated cases of finite abelian groups. Since inception , the study has been extended to some other important classes of finite abelian and nonabelian groups such as the dihedral , quaternion, semidihedral, and hamiltonian groups.
Other  different approaches have been so far, applied for the classification. The Fuzzy sets were introduced by Zadeh in 1965.  This theory of fuzzy sets has a wide range of applications, one of which is that of fuzzy groups developed by Rosenfield  in 1971. This by far, plays a pioneering role for the study of fuzzy algebraic structures. Other notions have been developed based on this theory. These, amongst others, include the notion of level subgroups by P.S. Das  used to characterize fuzzy subgroups of finite groups and that of equivalence of fuzzy subgroups introduced by Murali and Makamba  which we use in this work. ( Please, see [ 1 -  9 ]  )\\
By  the  way,  A group is nilpotent  if  it  has a normal  series  of  a finite length $ n.$

$$G = G_0 \geq G_1 \geq G_2 \geq \cdots \geq G_n = \{e\},   $$
where
$$G_i/G_{i+1} \leq Z(G/G_{i+1}).$$
By this notion, every finite $p$-group is nilpotent.  The nilpotence property is an hereditary one.  Thus,
\begin{enumerate}
\item [(i)] Any finite product of nilpotent group is nilpotent.
\item [(ii)] If $G$ is nilpotent of a class $c$, then, every subgroup and quotient group of $G$ is nilpotent and of class $\leq c$.
\end{enumerate}

The  problem of classifying the fuzzy subgroups of a finite group has so far experienced a very rapid progress. One particular case or the other have been treated by several papers such as the finite abelian as well as the non-abelian groups. The number of distinct fuzzy subgroups of a finite cyclic group of square-free order has been determined. Moreover,  a recurrence relation is indicated which can successfully be used to count the number of distinct fuzzy subgroups for two classes of finite abelian groups. They are  the arbitrary finite cyclic groups and finite elementary abelian $p$-groups. For the first class, the explicit formula obtained gave rise to an expression of a well-known central Delannoy numbers.
Some  forms of  propositions for classifying fuzzy subgroups for a class of finite $p$-groups have been made by Marius Tarnauceaus. It was from there, the study was extended to some important classes of finite non-abelian groups such as the dihedral and hamiltonian groups.  And thus, a method of determining the number and nature of fuzzy subgroups was developed with respect to the equivalence relation.
There are other different approaches  for the classification. The corresponding equivalence classes of fuzzy subgroups are closely connected to the chains of subgroups, and an essential role in solving counting problem is again played by the inclusion - exclusion principle.
This hereby leads to some recurrence relations, whose solutions have been easily found. For the purpose of using the Inclusion - Exclusion principle for generating the number of fuzzy subgroups, the finite $p$-groups has to be explored up to the maximal subgroups. The responsibility of describing the fuzzy subgroup structure of the finite nilpotent groups is the desired objective of this work.  Suppose that $(G, \cdot, e)$ is a group with identity $e$. Let $S(G)$ denote the collection of all fuzzy subsets of $G$. An element $\lambda \in S(G)$ is called  a fuzzy subgroup of $G$ whenever it satisfies some certain given  conditions . Such  conditions are  as  follows :

\begin{enumerate}
\item [(i)] $\lambda(ab) \geq \in \{\lambda(a), \lambda(b)\},\;\;\forall\; a,b\in G$;
\item [(ii)] $\lambda(a^{-1} \geq \lambda(a)$ for any $a \in G$.
\end{enumerate}
And, since $(a^{-1})^{-1}=a$, we have that $\lambda(a^{-1}) = \lambda(a)$, for any $a\in G$.\\
Also, by this notation and definition, $\lambda(e) = \sup\lambda(G)$.  [Marius [6]].\\
{\bf Theorem :}The set $FL(G)$ possessing all fuzzy subgroups of $G$ forms a lattice under the usual ordering of fuzzy set inclusion.  This is called the fuzzy subgroup lattice of $G$.

We define the level subset:
$$\lambda G_{\beta} = \{a\in G/\lambda(a) \geq \beta\}\;\;\mbox{for each $\beta\in [0,1]$}$$
The fuzzy subgroups of a finite $p$-group $G$ are thus, characterized, based on these subsets.  In the sequel, $\lambda$ is a fuzzy subgroup of $G$ if and only if its level subsets are subgroups in $G$.
This theorem gives a link between  $FL(G)$ and $L(G)$, the classical subgroup lattice of $G$.

Moreover, some natural relations on $S(G)$ can also be used in the process of classifying the fuzzy subgroups of a finite  $q$-group $G$.  One of them is defined by: $\lambda\sim \gamma$ iff $(\lambda(a)>\lambda(b) \Longleftrightarrow v(a)>v(b),\;\;\forall\;a,b\in G)$.  Alos, two fuzzy subgroups $\lambda, \gamma$ of $G$ and said to be distinct if $\lambda \times v$.

As a result of this development, let $G$ be a finite $p$-group and suppose that $\lambda :G \longrightarrow [0, 1]$ is a fuzzy subgroup of $G$.  Put $\lambda(G) = \{\beta_1,\beta_2,\dots, \beta_k\}$ with the assumption that $\beta_1 < \beta_2 > \cdots > \beta_k$.  Then, ends in $G$ is determined by $\lambda$.
$$\lambda G_{\beta_1} \subset \lambda G_{\beta_2} \subset \cdots \subset \lambda G_{\beta_k} = G \eqno(a)$$
Also, we have that:
$$\lambda(a) = \beta_t \Longleftrightarrow t = \max\{r/a \in \lambda G_{\beta_r}\}\Longleftrightarrow a\in \lambda G_{\beta_t} \backslash \lambda G_{\beta_{t-1}},\;$$
for any $a\in G$ and $t = 1,\dots, k$, where by  convention, set $\lambda G_{\beta_0}=\phi$.
\ \\ \ \\

 \section{  Methodology  } 

We are going to adopt a method  that will be used in counting the chains of fuzzy subgroups of an arbitrary finite $p$-group $G$ is described.
Suppose that $M_1, M_2, \dots, M_t$ are the maximal subgroups of $G$, and denote by $h(G)$ the number of chains of subgroups of $G$ which ends in $G$.  By simply applying the technique of computing $h(G)$, using the application of the Inclusion-Exclusion Principle, we have that:
$$h(G) = 2\left(\sum^t_{r=1}h(M_r) - \sum_{1\leq r_1<r_2\leq t}h(M_{r_1}\cap M_{r_2})\right.
\left.+\cdots+(-1)^{t-1}h\left(\bigcap^t_{r=1}M_r\right)\right)\hfill (\#)$$
In [6], (\#) was used to obtain the explicit formulas for some positive integers $n$.\\
{\bf Theorem  [ 1 ] [Marius]:}  The number of distinct fuzzy subgroups of a finite $p$-group of order $p^n$ which have a cyclic maximal subgroup
is:  \\   (i)     $h(\Z_{p^n}) = 2^n$,     (ii)    $\disp h(\Z_p\times \Z_{p^{n-1}}) = 2^{n-1}[2+(n-1)p]$

\ \\ \ \\
\bf { 3  \;\;	The  distinct   Number of The Fuzzy Subgroups  of The  Nilpotent Group of   $(D_{2^3} \times C_{2^m} ) $  For   $ m \geq  3 $ }   \\
\ \\ \ \\
{\bf    Proposition  1 ( see [ 13 ] ) :  } Suppose that $G =   \Z_4 \times  \Z_{2^{n}} , n \geq 2.$  Then,  $h(G) = 2^{n}[n^{2} + 5n - 2]$
\ \\ \ \\
{\bf    Proof :} $G$  has  three  maximal subgroups of which  two are isomorphic to   $ \Z_2 \times  \Z_{2^{n}}$  and  the   third is  isomorphic  to $ \Z_4 \times  \Z_{2^{n-1}}.$  \\
Hence,  $h( \Z_4 \times  \Z_{2^{n}})  = 2h( \Z_2 \times  \Z_{2^{n}}) + 2^{1}h( \Z_2 \times  \Z_{2^{n-1}}) +2^{2}h( \Z_2 \times  \Z_{2^{n-2}}) \\ + 2^{3}h( \Z_2 \times  \Z_{2^{n-3}})  + 2^{4}h( \Z_2 \times  \Z_{2^{n-4}})  + \cdots +  2^{n-2}h( \Z_2 \times  \Z_{2^2})  $\\
  
\begin{eqnarray*}
   = \;\;\;   2^{n+1}[2(n+1) + \sum^{n-2}_{j=1}[( n+1)-j ]    
\end{eqnarray*}                                        \\
=  $2^{n+1}[2(n+1) +\frac{1}{2}(n-2)(n+3)] =   2^{n}[n^{2} + 5n - 2] , n \geq 2$ \\
\ \\
We have that : $h( \Z_4 \times  \Z_{2^{n-1}}) =  2^{n-1}[(n-1)^{2} + 5(n-1) - 2]  \\  =  2^{n-1}[n^{2} + 3n - 6], n > 2  \hfill        \Box$
\ \\ \ \\
{\bf   Corrolary  1  :}     Following the last proposition, $ h( \Z_4 \times  \Z_{2^{5}}),h( \Z_4 \times  \Z_{2^{6}}),h( \Z_4 \times  \Z_{2^{7}}) $ and $h( \Z_4 \times  \Z_{2^{8}}) $ = 1536, 4096, 10496  and 26112  respectively. \\
\ \\ \ \\

{\bf   Theorem  A  ( see [ 15 ] ) :} Let $G =   D_{2^{n}} \times  \C_{2}$, the nilpotent group formed  by  the  cartesian product of  the  dihedral  group  of  order  $2^n$ and a cyclic  group  of  order  2. Then, the  number  of  distinct  fuzzy  subgroups  of  $G$  is  given  by : $h(G) = 2^{2n}(2n+1)- 2^{n+1} , n > 3  $ 
\ \\ \ \\ \ \\
{\bf Proof: } \\ The group $ D_{2^n}\times C_2,$   has one maximal subgroup which is isomorphic to  $\Z_2 \times \Z_{2^{n-2}}$, two  maximal subgroups which are isomorphic to $D_{2^{n-1}}\times C_2,$ and $2^2$ which are isomorphic to $D_{2^n}.$ \\ It thus, follows from the Inclusion-Exclusion Principle using equation, 
\begin{eqnarray*}\frac{1}{2}h(D_{2^n}\times C_2)  =  h(\Z_2\times \Z_{2^{n-1}})+4h(D_{2^n})-8h(D_{2^{n-1}})
-2h(\Z_2\times \Z_{2^{n-2}})+2h(D_{2^{n-1}}\times C_2)\end{eqnarray*}
By  recurrence relation principle  we have  :
 $$h(D_{2^n}\times C_2)  =  2^{2n}(2n+1)-2^{n+1},\;\;\;n>3$$   By the fundermental principle of mathematical induction, \\ set F(n) = $h(D_{2^n}\times C_2)$, assuming the truth of F(k) =$h(D_{2^k}\times C_2) = 2h(\Z_2 \times Z_{k-1})\\+ 8h(D_{2^k}-16hD_{2^{k-1}}-4h(\Z_2 \times \Z_{k-2})+ 4h(D_{2^{k-1}}\times C_2) = 2^{2k}(2k+1) - 2^{k+1},\\$ F(k+1)  = $h(D_{2^{k+1}}\times C_2) = 2h(\Z_2 \times \Z_{2^k}) + 8h(D_{2^{k+1}} - 16h(D_{2^k} - 4h(\Z_2 \times \Z_{k-1}) \\ + 4h(D_{2^k}\times C_2) = 2^{2}[2^{2k}(2k-3)-2^{k}],$ which is true.         $\hfill \Box$
\ \\ \ \\ \ \\
{\bf  Proposition  2  (  see [ 12 ] ) :  }  Suppose that $G =   D_{2^{n}} \times  \C_{4}.$ Then,  the  number of  distinct  fuzzy  subgroups  of  $G$  is  given  by : \begin{eqnarray*} 2^{2(n  - 2)}(64n  +  173)  +    3\sum^{n-3}_{j=1}2^{(n - 1 + j)}( 2n+1 -  2j   ) \end{eqnarray*} 
\ \\ \ \\
{\bf  Proof :}
\ \\ \ \\
$\frac{1}{2}h(D_{2^n} \times C_{4} ) = h(D_{2^n} \times C_{2} ) + 2h(D_{2^{n-1}} \times C_{4} ) - 4h(D_{2^{n-1}} \times C_{2} ) + h( \Z_4 \times  \Z_{2^{n-1}})   \\  -  2h( \Z_2 \times  \Z_{2^{n-1}}) - 2h( \Z_4 \times  \Z_{2^{n-2}})  + 8h( \Z_2 \times  \Z_{2^{n- 2}})  + h( \Z_{2^{n-1}}) - 4h( \Z_{2^{n-2}}) $ 
\ \\ \ \\
 $h(D_{2^n} \times C_{4} ) =   (n - 3).2^{2n + 2} +  2^{2(n  - 3)}(1460)  +  3[2^{n}(2n - 1) +2^{n + 1}(2n - 3) +2^{n + 2}(2n - 5)  + \cdots  + 7( 2^{2(n-2)})  ] $ \\
    \begin{eqnarray*}=  (n - 3).2^{2n + 2} +  2^{2(n  - 3)}(1460)  +
    3\sum^{n-3}_{j=1}2^{n - 1 + j}( 2n+1 -  2j   )   \\   =     2^{2(n  - 2)}(64n  +  173)  +    3\sum^{n-3}_{j=1}2^{n - 1 + j}( 2n+1 -  2j   )  
\end{eqnarray*}   
\ \\ \ \\ 
{\bf  Proposition  3  ( see  [10] ) :} Let $G$ be an abelian $p$-group of type $\Z_p\times\Z_p\times\Z_{p^n},$ where $p$ is a prime and $ n \geq 1.$ The number of distinct fuzzy subgroups of $G$ is \\
$h(\Z_p\times\Z_p\times\Z_{p^n}) = 2^{n}p(p+1)(n-1)(3+np+2p)
+(2^n-2)p^3-2^{n+1}(n-1)p^3+2^n[p^3+4(1+p+p^2)].$ \\
\ \\
{\bf Proof:}  There exist exactly $ 1 + p + p^2$ maximal subgroups for the abelian type  $\Z_p\times\Z_p\times\Z_{p^n},$ [Berkovich(2008)]. One of them is isomorphic to  \\ $\Z_p\times\Z_p\times\Z_{p^{n-1}},$ while each of the remaining $ p + p^2$ is isomorphic to $\Z_p\times\Z_{p^n}.$ Thus, by the application of the Inclusion-Exclusion Principle,we have as follows:
$h(\Z_p\times\Z_p\times\Z_{p^n}) 
 =   2^{n}p(p+1)(n-1)(3+np+2p)
+(2^n-2)p^3-2^{n+1}(n-1)p^3+2^n[p^3+4(1+p+p^2)]$  And thus,
\begin{eqnarray*}
 h(\Z_p\times\Z_p\times\Z_{p^{n-2}}) &=& 2^{n-2}[4+(3n-5)p+(n^2-5)p^2+(n^2-5n+8)p^3]-2p^2.
\end{eqnarray*}$\hfill \Box$   \\

{\bf   Corrolary  2  :} 
From (3)  above, obsreve that,  we  are  going  to  have  that:
\begin{eqnarray*}
\;\; h(\Z_3\times\Z_3\times\Z_{3^n})&=& 2^{n+1}[18n^2+9n+26]-54
\end{eqnarray*}
Similarly, for $p=5$, using the  same  analogy, we have
\begin{eqnarray*}
h(\Z_5\times\Z_5\times\Z_{5^n}) &=& 2[30h(\Z_5\times\Z_{5^n})+h(\Z_5\times\Z_5\times\Z_{5^{n-1}})\\
&&-p^3h(\Z_{5^n})-30h(\Z_{5^{n-1}})+125]
\end{eqnarray*}
And for $p=7$,

$h(\Z_7\times\Z_7\times\Z_{7^n}) = 2[56h(\Z_7\times\Z_{7^n})+h(\Z_7\times\Z_7\times\Z_{7^{n-1}})
-343h(\Z_{7^n})-56h(\Z_{7^{n-1}})+343]$

We have, in general,
$ h(\Z_p\times\Z_p\times\Z_{p^{n-2}}) =  2^{n-2}[4+(3n-5)p+(n^2-5)p^2+(n^2-5n+8)p^3]-2p^2$   $\hfill \Box$   \\
  \\



\  \\   \    \\ 
{\bf  Proposition  ( see  [14] ) :  } \\
Let   $     G =           (D_{2^3} \times C_{2^m} ) $  for   $ m \geq  3 $ .    Then ,  $    h( G)  =          m(89  - 23m)     +   (85)2^{m + 3}  -  124    $ \\
\ \\ \ 
{\bf    Proof  :  } \\
There  exist   seven  maximal subgroups , of which one is isomorphic to   $   D_{2^3} \times C_{2^m-1} , $  two being isomorphic to  $  C_{2^m} \times C_{2}   \times C_{2} ) ,  $ two isomorphic to  $  C_{2^m} \times C_{2}  , $ and  one each   isomorphic to $  C_{2^m} \times C_{4}  , $ and  $  C_{2^m}  $     respectively.  \\   Hence , by the inclusion  - exclusion principle,  using  the  propositions [1], [2],  [3], and Theorem  [1] we  have that : \\
\ \\ 
$\frac{1}{2}h(G) =  h(D_{2^3} \times C_{2^{m-1} })   +  2 h( C_{2^m} \times C_2 )  \times C_2 ) +  2 h( C_{2^m} \times C_2 ) + 2 h(C_{2^m}  \times C_4 ) +   h (C_{2^m} )  -   12 h (C_{2^m} \times C_2   ) -   6 h (C_{2^{m-1}}\times C_2   )\times C_2    )-   3 h (C_{2^{m-1}})  \times C_4   )  +   28 h (C_{2^{m-1}} \times C_2   )    +   2 h (C_{2^{m-1}}\times C_2   )\times C_2   )   +  4 h (C_{2^{m}}  \times C_2   )   +    h (C_{2^{m-1}} \times C_4   )  -   35 h (C_{2^{m-1}} \times C_2   ) -  7 h (C_{2^{m-1}} \times C_2   )   +   h(C_{2^{m-1}}  \times C_2   )     \\     =  h(D_{2^3} \times C_2^{m-1} )   +  2 h( C_{2^m} \times C_2 )  \times C_2 ) - 6 h( C_{2^m} \times C_2 ) + h(C_{2^m}  \times C_4 ) +   h (C_{2^m} )  -   4h (C_{2^{m-1}}\times C_2   )\times C_2    )  -   2 h (C_{2^{m-1}})  \times C_4   )  +   8 h (C_{2^{m-1}} \times C_2   )     $  \\  
  =   $   h(D_{2^3} \times C_2^{m-1} )    +  2^{m + 2 }(  6m^2   +  7m + 9)   - 32       -  (6)2^{m}(2m  +  2)  +      8m(2 ^m)  -  2^{m+2}{6m^2    -  5m  +  8}  +   2^6   +  2^{m}( m^2  +  5m  - 2 )       -    2^m  ( 3m  +  m^2  -    6  )   +   2^m    =   h(D_{2^3} \times C_2^{m-1} )    +  2^{m }(    46m -  4 )  +  2^m   +  32  =   h(D_{2^3} \times C_2^{m-1} )    +  2^{m  }(    46m -  3 )  +  32      $      \\
Hence , $ h(G) = 2 h(D_{2^3} \times C_2^{m-1} )    +  2^{m + 1 }(    46m -  3 )  +  64   $  = $   2^{m + 1 }(    46m -  3 )  +  64  + 2[ 2^{m  }(    46m -  49 )  +  64     +  2h(D_{2^3} \times C_2^{m-2} )    ]$     $  = 2^{m + 1 }(    46m -  3 )  +  64 +  2^{m + 1 }(    46m -  49 )  +  2^7 + 2^2 h(D_{2^3} \times C_2^{m-2} )     $    $   = 2^{m + 1 }(    46mm -  3 )  +  2^6  +  2^{m + 1 }(    46m -  49 )  +  2^7 + 2^2 [ 2^{m - 1 }(    46m -  95 )  +  64 + 2 h(D_{2^3} \times C_2^{m-3} )     $

 $h(D_{23n} \times C_{2^m} ) =   (46m - 3).2^{m + 1} +  2^6   +   (46m  - 49)2^{m + 1}  + 2^7 +    (46m  - 95)2^{m + 1}   +  2^8    +  2^3 h(D_{2^3} \times C_{2^{m - 3}} )  $ 
  \\    =  $  2^{m + 1}.[(46m - 3)  +  (46m  - 49)  +  (46m  -  95)]    +  2^6  +   2^7    +  2^8 +   2^3h(D_{2^3} \times C_2^{m-3} ) $    

    \[ \begin{array}{cc} =    &   \underbrace{    2^6   +   2^7   +  2^8  +   \cdots    +   2^{5+k}    }  \\ & \mbox { series (1) } \\    \mbox{  $ +  \;\;  2^{m + 1}.[46mk    +  $  } & \underbrace{    (  -  3   - 49  - 95   \cdots   ( -3  - 46(k - 1)))} ]    \\ & \mbox { series (2)  } \\   \mbox{  $    +  2^kh(D_{2^3} \times C_{2^m-k} ),   k  \in \{1 , 2 . 3.  \cdots  n \in  \ N  \}   $     }         \end{array} \] 

For the series  (1) , we have that,    $ U_m   =  2^6 . 2^{m-1}  =  2^{5+k}  ,  m + 5 =  k + 5 ,  \Rightarrow  m = k. $ We have that  $ S_{m = k } = 2^6[\frac{2^k  -  1}{2 - 1}] = 2^6 (2^k  .  1) $ \\

And  for the second series  (2), we have that ,  $ T_m  =  -  3 + (m-1)(-46) = -3 - 46(k-1)  \Rightarrow  m-1 =  k-1 ,  n=k  $  Hence ,  $ S_m=k  =  \frac{k}{2}[ 2(-3)  + (k-1)(-46) ]   =  \frac{k}{2} ( -6 - 46k + 46 )  = \frac{k}{2} ( 40 - 46k ) ,$  We have  that   
 $h(D_{23n} \times C_{2^m} ) = \frac{k}{2} ( 40 - 46k ) +  2^6 (2^k  .  1 ) +  2^k h(D_{3} \times C_{2^m-k} $.  By setting   $ m = k $ we have  that  $ k = m - 3$.Hence ,  $ h(D_{2^3} \times C_{2^m} ) =  (m - 3)  (20  -  23m) + 2^6  (2^{m - 3}  -  1)  +    2^m-3  h(D_{3} \times C_{2^3}) $

$h(G) =    (m - 3)(20  -  23m) + 2^6  (2^{m - 3}  -  1)    +  2^{m - 3}(5376)  $   =   $ (m - 3)(20  -  23m) + 2^{m - 3}  -  2^6      +  2^{m + 5}(21)  $     = $   {20m  - 23m^2  -   60  +  69m}  +  2^{m + 3} -  2^6  +   (21)2^{m + 5}  $  = $ (89m  - 23m^2  -   60 )  +  2^{m + 3} -  2^6   +   (21)2^{m + 5}  $        =  $      m(89  - 23m)  -  124    +   (85)2^{m + 3}       \hfill      \Box         $   \\    


\ \\ \ \\

{\bf    Theorem  ( see  [11] ) :}  Let  $G =   \Z_{2^{n}}\times \Z_8$, then  $ h(G)  =    \frac{1}{3}(2^{n+1})(n^{3} + 12n^{2} + 17n  -  24 )$ \\
{\bf    Proof :} The   three  maximal subgroups of $G$ have the following properties :  \\
one is isomorphic to $\Z_8 \times  \Z_{2^{n-1}})$, while two are isomorphic to  $\Z_4 \times  \Z_{2^n})$ .\\
 We have :  $ \frac{1}{2}h(G) = 2h( \Z_4 \times  \Z_{2^n})  + h( \Z_8 \times  \Z_{2^{n-1}})  - 3 h( \Z_4 \times  \Z_{2^{n-1}}) + h( \Z_4 \times  \Z_{2^{n-1}})\\ = 2h( \Z_4 \times  \Z_{2^n})  + h( \Z_8 \times  \Z_{2^{n-1}})  - 2 h( \Z_4 \times  \Z_{2^{n-1}}) \\ = h( \Z_8 \times  \Z_{2^{n-1}})  + 2 h( \Z_4 \times  \Z_{2^n})  -  h( \Z_4 \times  \Z_{2^{n-1}})$  \\
Hence ,  $  h(G)= 4h( \Z_4 \times  \Z_{2^n})   - 4 h( \Z_4 \times  \Z_{2^{n-1}})   + 2 h( \Z_8 \times  \Z_{2^{n-1}})\\
 = 4h( \Z_4 \times  \Z_{2^n})  +4 h( \Z_4 \times  \Z_{2^{n-1}})   + 8 h( \Z_4 \times  \Z_{2^{n-2}}) - 16 h( \Z_4 \times  \Z_{2^{n-3}}) \\ + 32h( \Z_4 \times  \Z_{2^{n-3}})  - 32 h( \Z_4 \times  \Z_{2^{n-4}})  + 16 h( \Z_8 \times  \Z_{2^{n-4}})  \\ =  4h( \Z_4 \times  \Z_{2^n})  + 4h( \Z_4 \times  \Z_{2^{n-1}})  + 8h( \Z_4 \times  \Z_{2^{n-2}}) + 16h( \Z_4 \times  \Z_{2^{n-3}}) \\ + 32h( \Z_4 \times  \Z_{2^{n-4}})  - 64h( \Z_4 \times  \Z_{2^{n-5}})   + 32 h( \Z_8 \times  \Z_{2^{n-5}})  +  \cdots    - 2^{j+1}h( \Z_4 \times  \Z_{2^{n-j}})  \\  + 2^{j}h( \Z_8 \times  \Z_{2^{n-j}})  $ , for  $ n - j = 3 $     \begin{eqnarray*}  = 4h( \Z_4 \times  \Z_{2^n})  + 2^{n-3}h( \Z_8 \times  \Z_{2^3}) -  2^{n-1}h( \Z_4 \times  \Z_{2^3})  +   \sum^{n-3}_{k=1}[ 2^{k+1}h( \Z_4 \times  \Z_{2^{n-k}})   \end{eqnarray*}  
$ =  2^{n+2}[ n^{2} + 5n + 3  ] +  \sum^{n-3}_{k=1}h( \Z_4 \times  \Z_{2^{n-k}}) $
 = $  2^{n+2}((n^{2}+ 5n+ 3) + \frac{1}{6}(n-3)(n^{2} + 9n + 14 ))  \\ = \frac{1}{3}(2^{n+1})(n^{3} + 12n^{2} + 17n  -  24 ) ,   n > 2  .   \hfill   \Box$  \\

{\bf  Proposition   ( see  [16] :  }  Suppose that $G =   D_{2^{n}} \times  \C_{8}.$ Then,  the  number of  distinct  fuzzy  subgroups  of  $G$  is  given  by : \begin{eqnarray*}   2^{2(n-1)}(  6n + 113)  +     2^{n}[ 13 -  6n - 2n^{2}  +   3\sum^{n-3}_{j=1}2^{( j - 1j)}( 2n+1 -  2j )]  \\  +  \frac{1}{3}(2^{n+2})[(n-1)^{3} + (n-2)^{3} + 24n^{2} - 38n - 30 +  \sum^{n-5}_{k=1}2^{k}[ ( n -  2 - k )^{3} + 12(n-2-k)^{2} + 17(n-k) - 58 ]]  \end{eqnarray*} 
{\bf  Proof :}  $ h(D_{2^n} \times C_{8} ) = 2h(\Z_{2^{n-1}}) + 2h(D_{2^n} \times Z_{4} ) + 2h(D_{2^{n-1}} \times C_{8} )   \\ + 4h(Z_{2^{n-2}} \times C_{8} ) + 2^{4}h(Z_{2^{n-3}} \times C_{8} ) + 2^{6}h(\Z_{2^{n-4}} \times C_{8} ) - 2^{8}h( \Z_{2^{n-5}} \times  \Z_{2^{3}})  \\   - 4h(  \Z_{2^{n-1}} \times  \Z_{2^{2}})  + 2^{10}h( \Z_{2^{n-5}}) \times \Z_{2^{2}} - 2^{9}h( \Z_{2^{n-5}}) -  2^{9}h(D_{2^{n-4}} \times C_{2^{2}} )  \\  + 2^{8}h(D_{2^{n-4}} \times C_{2^{3}} )  \\ =    2^{n}  +  2h(D_{2^n} \times C_{4} ) + 2h( \Z_{2^{n-1}} \times  \Z_{2^{3}})  + 2^{2}h( \Z_{2^{n-2}} \times  \Z_{2^{3}})  \\  - 2^{2(n-3)}h( \Z_{2^{2}} \times  \Z_{2^{3}})  +  2^{2(n-2)}h( \Z_{2^{2}} \times  \Z_{2^{2}}  - 2^{2}h( \Z_{2^{n-1}} \times  \Z_{2^{2}})  - 2^{2n-5}h( \Z_{2^{2}}) \\  -  2^{2n-5}h( D_{2^{3}} \times  \Z_{2^{2}})   +  2^{2(n-3)}h( D_{2^{3}} \times  \Z_{2^{3}})   $ \\
    \begin{eqnarray*} + \;\;\; \;\;     3\sum^{n-5}_{i=1}2^{2ij}h( \Z_{2^{n-2-i}} \times  \Z_{2^{3}})  
\end{eqnarray*}     as  required.  $    \hfill        \Box $

{\bf  Theorem  : }Let  $ G = D_{2^{4}} \times \C_{2^{4}} .$ Then , $ h(  G ) =  61384 $ \\

{\bf  Proof   : } There exist seven maximal subgroups .  Two  isomorphic to  $ D_{2^{4}} \times \C_{2^{3}} .$ two  isomorphic to  $ D_{2^{3}} \times \C_{2^{4}} .$ two  isomorphic to  $ D_{2^{4}} \times \C_{2^{2}} ,$ while the seventh is isomorphic to $  \Z_{2^{4}} .$ \\

Hence , we have that : $\frac{1}{2}h(  G ) =  2h(D_{2^4} \times Z_{2^2} ) + 2h(D_{2^{4}} \times Z_{2^3} ) + 2h(D_{2^{3}} \times Z_{2^4} ) -  6h( D_{2^3} \times \Z_{2^{3}})  - 6h( \Z_{2^4} \times \Z_{2^{2}}) - 3h( \Z_{2^{3 }} \times \Z_{2^3})  -6h( \Z_{2^4}) + 2h( D_{2^{3}} \times \Z_{2^{3}})  + 28h(\Z_{2^3} \times Z_{2^2} ) + 2h(Z_{2^{4}} \times Z_{2^2} ) + 2h( \Z_{2^4} ) + h( Z_{2^3} \times \Z_{2^{3}})  - 35h( \Z_{2^3} \times \Z_{2^{2}})  + 21h( \Z_{2^{3 }} \times \Z_{2^2})  -7h( \Z_{2^3} \Z_{2^2}) + h(\Z_{2^3} \times  \Z_{2^2}) $ \\  
= $ 2[ h(D_{2^4} \times Z_{2^2} ) + h(D_{2^{4}} \times Z_{2^3} ) + h(D_{2^{3}} \times Z_{2^4} ) -  2h( D_{2^3} \times \Z_{2^{3}})  - 2h( \Z_{2^4} \times \Z_{2^{2}}) - h( \Z_{2^{3 }} \times \Z_{2^3})  + 4h( D_{2^{3}} \times \Z_{2^{2}})  -3h(\Z_{2^4} ) +  \frac{1}{2}h(Z_{2^{4}} ) ] $  \\

$ \therefore    h(G) = 4[ 700 + 8416 + 10744 -10752 – 1088 + 162 + 704 – 40 ]     \\    = 4[ 15346 ] = 61384  \hfill  \Box $  \\

\section{  Computation   for    $G = D_{2^{4}} \times \C_{2^{n}} ,  n \geq  4.$ }
Our  computation  on  the  algebraic  fuzzy   structure  given   actually   has  an  outcome  which  involves  multiple  sums  \\
\ \\ 
{\bf Proof :}
\ \\ 
The  maximal  subgroups  are  :    \\    $ (D_{2^4} \times C_{2^{n -1}} ) , 2(D_{2^3} \times C_{2^{n}} ) , 2(D_{2^n} \times C_{2^2} ),  (D_{2^n} \times C_{2^3} )  $ and  $ ( C_{2^n} ) . $ \\
We  have  that : 
$\frac{1}{2}h( G ) = h(D_{2^4} \times C_{n-1} )  +   2h(D_{2^3} \times C_{n} ) + 2h(D_{2^{n}} \times C_{2^2} ) + h(D_{2^{n}} \times C_{2^3} )  + h( C_{2^n} )    - 6h( D_{2^3} \times \Z_{2^{n-1}}) - 6h( \Z_{2^n} \times \Z_{2^{2}}) - 3h( \Z_{2^{n – 1 }} \times \Z_{2^3})  - 6h(  \Z_{2^n})  + 2h(D_{2^3} \times C_{2^{n-1}} )  + 28h(C_{2^{n-1}} \times C_{2^{n}} )   +  h(C_{2^{n-1}} \times C_{2^{3}} )  +  2h(C_{2^n} \times C_{2^2} )   + 2h( \Z_{2^{n}})  - 35h(C_{2^{n-1}} \times C_{2^2} )   +  21h(C_{2^{n-1}} \times C_{2^2} )    -  7h(C_{2^{n-1}} \times C_{2^2} )   +  h(C_{2^{n-1}} \times C_{2^2} )  $ \\

$ = h(D_{2^{4}} \times C_{2^{n-1}} )  +  2h(D_{2^3} \times C_2^{n} ) + 2h(D_{2^{n}} \times C_{2^2} ) + h(D_{2^{n}} \times C_{2^3} ) - 4h( D_{2^3} \times \Z_{2^{n-1}}) - 4h( \Z_{2^n} \times \Z_{2^{2}}) - 2h( \Z_{2^{n – 1 }} \times \Z_{2^3}) + 8h( \Z_{2^{n-1}} \times \Z_{2^2})  -  3h( \Z_{2^{n}}) $ \\

\begin{displaymath}    \frac{1}{2}h(  G ) =  h( D_{2^4} \times \Z_{2^{n-k}} )  +  2h( D_{2^3} \times \Z_{2^n} )   - 4h( D_{2^3} \times \Z_{2^{n-k}})  -  4h( \Z_{2^n} \times \Z_{2^2}) \end{displaymath}    \begin{displaymath}  -  2h( \Z_{2^{n-k}} \times \Z_{2^3})    + 8h( \Z_{2^{n-k}} \times \Z_{2^2})   +    \sum_{j = 1}^{k} h( D_{2^{n-1+j}} \times \Z_{2^3}) + 2\sum_{j = 1}^{k} h( D_{2^{n-1+j}} \times \Z_{2^2}) - 3\sum_{j = 1}^{k} h( \Z_{2^{n+1-j}}) \end{displaymath} \begin{displaymath} - 2\sum_{j = 1}^{k-1} h( D_{2^3} \times \Z_{2^{n-j}}) + 4\sum_{j = 1}^{k-1} h( D_{2^{n-j}} \times \Z_{2^2}) - 2\sum_{j = 1}^{k-1} h( D_{2^n-j} \times \Z_{2^3}), \end{displaymath}
 whence , $ n -  k = 4 ,  \Rightarrow  k  = n -  4.   \therefore  h(G)  =  2h( D_{2^4} \times \Z_{2^4}) +  4h( D_{2^3} \times \Z_{2^n})   - 8h( D_{2^3} \times \Z_{2^4}) -  8h( \Z_{2^n} \times \Z_{2^n}) -  4h( \Z_{2^4} \times \Z_{2^3})  + 16h( \Z_{2^4} \times \Z_{2^2}) +  $

\begin{displaymath}   2\sum_{j = 1}^{n-4} h( D_{2^{n-1+j}} \times \Z_{2^3})    +  4\sum_{j = 1}^{n-4} h( D_{2^{n-1+j}} \times \Z_{2^2})  -  6\sum_{j = 1}^{n-4} h(  \Z_{2^{n+1-j}})   \end{displaymath}  \begin{displaymath} -  4\sum_{j = 1}^{n-5} h( D_{2^3} \times \Z_{2^{n-j}})   +  8\sum_{j = 1}^{n-5} h( D_{2^n-j} \times \Z_{2^2})   -  4\sum_{j = 1}^{n-5} h( D_{2^{n-j}} \times \Z_{2^3})    \end{displaymath}

\begin{displaymath} \therefore  h(G)  =  2^{n+3}(422  -  n^2    -  5n    )   -  9n^2    +    356n   -   29160  +     2\sum_{j = 1}^{n-4} h( D_{2^{n-1+j}} \times \Z_{2^3})  \end{displaymath}  \begin{displaymath}  +  4\sum_{j = 1}^{n-4} h( D_{2^{n-1+j}} \times \Z_{2^2})  -  6\sum_{j = 1}^{n-4} h(  \Z_{2^{n+1-j}})    -  4\sum_{j = 1}^{n-5} h( D_{2^3} \times \Z_{2^{n-j}})   +  8\sum_{j = 1}^{n-5} h( D_{2^n-j} \times \Z_{2^2})   -  4\sum_{j = 1}^{n-5} h( D_{2^{n-j}} \times \Z_{2^3})    \end{displaymath}

\begin{displaymath}   =  2^{n+3}(422  -  n^2    -  5n    )   -  9n^2    +    356n   -   29160  +     \sum_{j = 1}^{n-4}[2h( D_{2^{n-1+j}} \times \Z_{2^3})    +  4h( D_{2^{n-1+j}} \times \Z_{2^2})  -  6h(  \Z_{2^{n+1-j}})   ]\end{displaymath}  \begin{displaymath}  -  \sum_{j = 1}^{n-5}[4 h( D_{2^3} \times \Z_{2^{n-j}})   -  8h( D_{2^n-j} \times \Z_{2^2})   +  4h( D_{2^{n-j}} \times \Z_{2^3}) ]   \end{displaymath}

\ \\ \

Hence , proved  as  required  $  \hfill    \Box    $

 \section{	Applications  }

        The  computations  so far  by the use of   GAP ( General  AlgorithmAlgorithms  and  Programming  )  and  the  Inclusion  - Exclusion  Principle  can be certified here as being very useful in the  computations  of the distinct  number of fuzzy subgroups for the finite  nilpotent   $ p $ - groups .

\ \\ \ \\

 \section{	INSTANCES   }

We have the following examples  as  parts  surfacing from our computations  so  far. The  readers  may  consider  the  examples below  in  tabular  format. \\
\ \\ \ 
{\bf   \;\;   Example  1 :}   \\
\section*{Table $1$\\
Table Summarizing     some  Number of    Distinct Fuzzy Subgroups of  $(D_{2^3} \times C_{2^n} ) $  FOR  $  \geq  3 $ }
\begin{center}
\begin{tabular}{|c|c|c|c|c|c|c|c|c |c|c|c |c| |c| |c|c  |c| |c|c|c| c|c|c| |c| |c|  }                                                                                                                                                                                                                    \hline
S/N  for  the    Number of  $m$ & \multicolumn{1}{|c|}{3} & \multicolumn{1}{|c|}{4} &     \multicolumn{1}{|c|}{5}&\multicolumn{1}{|c|}{6}&\multicolumn{1}{|c|}{7}&\multicolumn{1}{|c|}{8}&\multicolumn{1}{|c|}{9}&\multicolumn{1}{|c|}{10}  \\ \hline
$h(G)     =     (D_{2^3} \times C_{2^n} ) ,  n \geq  3  $ & 5376 & 10728 & 21506 & 43347 & 86536 & 173320 & 347098 & 694910        \\ \hline
\end{tabular}
\end{center}
\ \\ \ \\ \ \\
{\bf   \;\;   Example  2 :}  Now,  since  the stipulated   condition  that  $ m  \geq  3$  must   definitely  be  fulfilled then  the  readers  may  consider  the  examples below  in  tabular  format.  \\
\section*{Table $2$\\
Table Summarizing   some  Number of    Distinct Fuzzy Subgroups of  $(D_{2^4} \times C_{2^n} ) $  FOR  $ n \geq  4 $ }
\begin{center}
\begin{tabular}{|c|c|c|c|c|c|c|c|c |c|c|c |c| |c| |c|c |c| |c|c|c| c|c|c| |c| |c| } \hline
S/N for the Number of $n$ & \multicolumn{1}{|c|}{4} & \multicolumn{1}{|c|}{5} & \multicolumn{1}{|c|}{6}  \\ \hline
$h(G) = (D_{2^4} \times C_{2^n} ),   n \geq  4   $ & 20, 200& 375, 648& 3, 893, 800 \\ \hline

\end{tabular}
\end{center}

 \section{		Conclusion }
\ \\
The    discoveries  from  our  studies  so far ,  has  helped to   observe  that any finite product of nilpotent group is nilpotent.  Also, the problem of classifying the fuzzy subgroups of a finite group has  experienced a very rapid progress.     Finally,  the method can be used in further computations  up to the generalizations  of similar and other given structures  
\\\\ \\ \
\bf { Funding: {   This research received no external funding
\  \\ \ \\



{\bf    Competing of interests statement : } The authors declare that in this paper,  there is no competing of interests.  \\


\  \\  \  \\  \  \\  \ 
{\bf References}
\begin{description}
\item [(1)]   M.  Mashinchi   and    M.  Mukaidono (1992). A classification of fuzzy subgroups. Ninth Fuzzy  System Symposium, Sapporo, Japan, 649-652.
\item [(2)]   M.    Mashinchi  ,  M.    Mukaidono   .(1993). On fuzzy subgroups  classification, Research Reports of  Meiji Univ. (9) , 31-36.
\item [(3)]  V.  Murali   and   B.  B.  Makamba    (2003).  On an equivalence of Fuzzy Subgroups III, Int. J. Math. Sci. 36, 2303-2313.
\item [(4)] Odilo Ndiweni (2014).  The classification of Fuzzy subgroups of the Dihedral Group $D_n$, for $n$, a product of distinct primes.  A Ph.D. thesis, Univ. of Fort Hare, Alice, S.A.
\item [(5)]  M.  Tarnauceanu   (2009). The number of fuzzy subgroups of finite cyclic groups and Delannoy numbers, European J. Combin. (30), 283-289, doi: 10.1016/j.ejc.2007.12.005.
\item [(6)] M.  Tarnauceanu    (2011).  Classifying fuzzy subgroups for a class of finite $p$-groups.  ``ALL CUZa'' Univ. Iasi, Romania.
\item [(7)] M.  Tarnauceanu   (2012). Classifying  fuzzy subgroups of finite nonabelian groups. Iran.J.Fussy Systems. (9) 33-43
\item [(8)] M. Tarnauceanu ,   L.  Bentea   (2008). A note on  the number of  fuzzy subgroups of finite groups, Sci. An. Univ. "ALL.Cuza" Iasi, Matt., (54) 209-220.
\item [(9)]  M. Tarnauceanu,   L.  Bentea  (2008). On  the number of fuzzy subgroups of finite abelian groups, Fuzzy Sets  and Systems (159), 1084-1096, doi:10.1016/j.fss.2017.11.014 
\item [(10)] 	S. A. ADEBISI    and M. EniOluwafe (2020) The Abelian Subgroup :    $ \Z_p \times \Z_p \times  \Z_{p^n} $    ,  $p$ is Prime and  $ n \geq  1 $  .  Progress in Nonlinear Dynamics and Chaos Vol. 7, No. 1, 2019, 43-45 ISSN: 2321 – 9238 (online) Published on 21 September 2019   www.researchmathsci.orgDOI: http://dx.doi.org/10.22457/pindac.80v7n1a343
\item [(11)] S. A. ADEBISI, M. OGIUGO AND M. ENIOLUWAFE (2020) THE FUZZY SUBGROUPS FOR THE ABELIAN    STRUCTURE: $ h( \Z_{8} \times \Z_{2^n}) $ for  $ n > 2.$  Journal of the Nigerian Mathematical Society, Vol. 39, Issue 2, pp. 167-171.

\item [(12)] S. A. ADEBISI, M. Ogiugo  and M. EniOluwafe(2020)Computing the Number of Distinct Fuzzy Subgroups for the Nilpotent $ p$-Group of $  D_{2^n}\times C_4 $  International J.Math.Combin.1(2020),86-89.
\item [(13)] 	S. A. ADEBISI, M. Ogiugo  and M. EniOluwafe  (2020) Determining  The Number Of Distinct Fuzzy  Subgroups  For  The  Abelian Structure    :  $ \Z_4 \times \Z_{2^{n-1}} $   ,  $n >2 $.Transactions of the Nigerian Association of Mathematical Physics Volume 11, (January - June, 2020 Issue), pp 5 - 6
\item [(14)] S. A. ADEBISI ,   M. OGIUGO and  M. ENIOLUWAFE (2022)THE FUZZY SUBGROUPS FOR THE  NILPOTENT  (  $p $-GROUP) OF $ (D_{2^3} \times  C_{2^m} )  $ FOR  $m  \geq 3 $  Journal of Fuzzy Extension and Applications    www.journal-fea.com  ( Accepted  for  publication )
\item [(15)] S. A. ADEBISI  and  M. EniOluwafe (2020) An explicit formula for the number of distinct fuzzy subgroups of the Cartesian product of the dihedral group of order 2n with a cyclic group of order 2 Universal Journal of Mathematics and Mathematical Sciences.http://www.pphmj.com http://dx.doi.org/10.17654/UM013010001 Volume 13, no1, 2020, Pages 1-7 ISSN: 2277-1417. \\  (http://www.pphmj.com/journals/articles/1931.htm)

\item [(16)]    S. A. ADEBISI, M. Ogiugo  and M. EniOluwafe(2020) On the $ p$-Groups of the Algebraic Structure of $  D_{2^n} \times  C_8 $ International  J.Math. Combin. Vol.3(2020), 100-103

\end{description}
\end{document}